\documentclass[12pt,a4paper]{article}

\usepackage{amsfonts,amssymb,amsthm,amsmath,latexsym}
\usepackage{subeqnarray,indent,url,cite,bm}

\date{October 31, 2019 \\[1mm] revised January 6, 2020}

\begin{document}

\title{A remark on the enumeration of rooted labeled trees}

\author{
     {\small Alan D.~Sokal}                  \\[2mm]
     {\small\it Department of Mathematics}   \\[-2mm]
     {\small\it University College London}   \\[-2mm]
     {\small\it Gower Street}                \\[-2mm]
     {\small\it London WC1E 6BT}             \\[-2mm]
     {\small\it UNITED KINGDOM}              \\[-2mm]
     {\small\tt sokal@math.ucl.ac.uk}        \\[-2mm]
     {\protect\makebox[5in]{\quad}}  
     \\[-2mm]
     {\small\it Department of Physics}       \\[-2mm]
     {\small\it New York University}         \\[-2mm]
     {\small\it 726 Broadway}                \\[-2mm]
     {\small\it New York, NY 10003}          \\[-2mm]
     {\small\it USA}                         \\[-2mm]
     {\small\tt sokal@nyu.edu}               \\[3mm]
}

\maketitle
\thispagestyle{empty}   

\begin{abstract}
Two decades ago, Chauve, Dulucq and Guibert showed that
the number of rooted trees on the vertex set $[n+1]$
in which exactly $k$ children of the root are lower-numbered than the root
is $\binom{n}{k} \, n^{n-k}$.
Here I give a simpler proof of this result.
\end{abstract}


\bigskip
\bigskip
\noindent
{\bf Mathematics Subject Classification (MSC 2010) codes:}
05A15 (Primary);  05A10, 05A19, 05C05 (Secondary).

\clearpage

\newtheorem{theorem}{Theorem}[section]
\newtheorem{proposition}[theorem]{Proposition}
\newtheorem{lemma}[theorem]{Lemma}
\newtheorem{corollary}[theorem]{Corollary}
\newtheorem{definition}[theorem]{Definition}
\newtheorem{conjecture}[theorem]{Conjecture}
\newtheorem{question}[theorem]{Question}
\newtheorem{problem}[theorem]{Problem}
\newtheorem{example}[theorem]{Example}

\renewcommand{\theenumi}{\alph{enumi}}
\renewcommand{\labelenumi}{(\theenumi)}
\def\eop{\hbox{\kern1pt\vrule height6pt width4pt
depth1pt\kern1pt}\medskip}
\def\prf{\par\noindent{\bf Proof.\enspace}\rm}
\def\rmk{\par\medskip\noindent{\bf Remark\enspace}\rm}

\newcommand{\textbfit}[1]{\textbf{\textit{#1}}}

\newcommand{\bigdash}{%
\smallskip\begin{center} \rule{5cm}{0.1mm} \end{center}\smallskip}

\newcommand{\safepar}{ {\protect\hfill\protect\break\hspace*{5mm}} }

\newcommand{\be}{\begin{equation}}
\newcommand{\ee}{\end{equation}}
\newcommand{\<}{\langle}
\renewcommand{\>}{\rangle}
\newcommand{\widebar}{\overline}
\def\reff#1{(\protect\ref{#1})}
\def\spose#1{\hbox to 0pt{#1\hss}}
\def\ltapprox{\mathrel{\spose{\lower 3pt\hbox{$\mathchar"218$}}
    \raise 2.0pt\hbox{$\mathchar"13C$}}}
\def\gtapprox{\mathrel{\spose{\lower 3pt\hbox{$\mathchar"218$}}
    \raise 2.0pt\hbox{$\mathchar"13E$}}}
\def\textprime{${}^\prime$}
\def\proof{\par\medskip\noindent{\sc Proof.\ }}
\def\firstproof{\par\medskip\noindent{\sc First Proof.\ }}
\def\secondproof{\par\medskip\noindent{\sc Second Proof.\ }}
\def\alternateproof{\par\medskip\noindent{\sc Alternate Proof.\ }}
\def\algebraicproof{\par\medskip\noindent{\sc Algebraic Proof.\ }}
\def\combinatorialproof{\par\medskip\noindent{\sc Combinatorial Proof.\ }}
\def\proofof#1{\bigskip\noindent{\sc Proof of #1.\ }}
\def\firstproofof#1{\bigskip\noindent{\sc First Proof of #1.\ }}
\def\secondproofof#1{\bigskip\noindent{\sc Second Proof of #1.\ }}
\def\thirdproofof#1{\bigskip\noindent{\sc Third Proof of #1.\ }}
\def\algebraicproofof#1{\bigskip\noindent{\sc Algebraic Proof of #1.\ }}
\def\combinatorialproofof#1{\bigskip\noindent{\sc Combinatorial Proof of #1.\ }}
\def\sketchofproof{\par\medskip\noindent{\sc Sketch of proof.\ }}
\renewcommand{\qed}{ $\square$ \bigskip}
\newcommand{\myendremark}{ $\blacksquare$ \bigskip}
\def\half{ {1 \over 2} }
\def\third{ {1 \over 3} }
\def\twothird{ {2 \over 3} }
\def\smfrac#1#2{{\textstyle{#1\over #2}}}
\def\smhalf{ {\smfrac{1}{2}} }
\newcommand{\real}{\mathop{\rm Re}\nolimits}
\renewcommand{\Re}{\mathop{\rm Re}\nolimits}
\newcommand{\imag}{\mathop{\rm Im}\nolimits}
\renewcommand{\Im}{\mathop{\rm Im}\nolimits}
\newcommand{\sgn}{\mathop{\rm sgn}\nolimits}
\newcommand{\tr}{\mathop{\rm tr}\nolimits}
\newcommand{\supp}{\mathop{\rm supp}\nolimits}
\newcommand{\disc}{\mathop{\rm disc}\nolimits}
\newcommand{\diag}{\mathop{\rm diag}\nolimits}
\newcommand{\tridiag}{\mathop{\rm tridiag}\nolimits}
\newcommand{\AZ}{\mathop{\rm AZ}\nolimits}
\newcommand{\perm}{\mathop{\rm perm}\nolimits}
\def\hboxscript#1{ {\hbox{\scriptsize\em #1}} }
\renewcommand{\emptyset}{\varnothing}
\newcommand{\eqdef}{\stackrel{\rm def}{=}}

\newcommand{\restrict}{\upharpoonright}

\newcommand{\compinv}{{\langle -1 \rangle}}   

\newcommand{\scra}{{\mathcal{A}}}
\newcommand{\scrb}{{\mathcal{B}}}
\newcommand{\scrc}{{\mathcal{C}}}
\newcommand{\scrd}{{\mathcal{D}}}
\newcommand{\scre}{{\mathcal{E}}}
\newcommand{\scrf}{{\mathcal{F}}}
\newcommand{\scrg}{{\mathcal{G}}}
\newcommand{\scrh}{{\mathcal{H}}}
\newcommand{\scri}{{\mathcal{I}}}
\newcommand{\scrj}{{\mathcal{J}}}
\newcommand{\scrk}{{\mathcal{K}}}
\newcommand{\scrl}{{\mathcal{L}}}
\newcommand{\scrm}{{\mathcal{M}}}
\newcommand{\scrn}{{\mathcal{N}}}
\newcommand{\scro}{{\mathcal{O}}}
\newcommand\scroo{
  \mathchoice
    {{\scriptstyle\mathcal{O}}}
    {{\scriptstyle\mathcal{O}}}
    {{\scriptscriptstyle\mathcal{O}}}
    {\scalebox{0.6}{$\scriptscriptstyle\mathcal{O}$}}
  }
\newcommand{\scrp}{{\mathcal{P}}}
\newcommand{\scrq}{{\mathcal{Q}}}
\newcommand{\scrr}{{\mathcal{R}}}
\newcommand{\scrs}{{\mathcal{S}}}
\newcommand{\scrt}{{\mathcal{T}}}
\newcommand{\scrv}{{\mathcal{V}}}
\newcommand{\scrw}{{\mathcal{W}}}
\newcommand{\scrz}{{\mathcal{Z}}}

\newcommand{\bfa}{{\mathbf{a}}}
\newcommand{\bfb}{{\mathbf{b}}}
\newcommand{\bfc}{{\mathbf{c}}}
\newcommand{\bfd}{{\mathbf{d}}}
\newcommand{\bfe}{{\mathbf{e}}}
\newcommand{\bfh}{{\mathbf{h}}}
\newcommand{\bfj}{{\mathbf{j}}}
\newcommand{\bfi}{{\mathbf{i}}}
\newcommand{\bfk}{{\mathbf{k}}}
\newcommand{\bfl}{{\mathbf{l}}}
\newcommand{\bfm}{{\mathbf{m}}}
\newcommand{\bfn}{{\mathbf{n}}}
\newcommand{\bfp}{{\mathbf{p}}}
\newcommand{\bfx}{{\mathbf{x}}}
\newcommand{\bfy}{{\mathbf{y}}}
\renewcommand{\k}{{\mathbf{k}}}
\newcommand{\n}{{\mathbf{n}}}
\newcommand{\vv}{{\mathbf{v}}}
\newcommand{\bv}{{\mathbf{v}}}
\newcommand{\w}{{\mathbf{w}}}
\newcommand{\x}{{\mathbf{x}}}
\newcommand{\y}{{\mathbf{y}}}
\newcommand{\cc}{{\mathbf{c}}}
\newcommand{\zero}{{\mathbf{0}}}
\newcommand{\one}{{\mathbf{1}}}
\newcommand{\bmm}{{\mathbf{m}}}

\newcommand{\ahat}{{\widehat{a}}}
\newcommand{\Zhat}{{\widehat{Z}}}

\newcommand{\C}{{\mathbb C}}
\newcommand{\D}{{\mathbb D}}
\newcommand{\Z}{{\mathbb Z}}
\newcommand{\N}{{\mathbb N}}
\newcommand{\Q}{{\mathbb Q}}
\newcommand{\PP}{{\mathbb P}}
\newcommand{\R}{{\mathbb R}}
\newcommand{\RR}{{\mathbb R}}
\newcommand{\E}{{\mathbb E}}

\newcommand{\Sym}{{\mathfrak{S}}}
\newcommand{\SymB}{{\mathfrak{B}}}
\newcommand{\Alt}{{\mathrm{Alt}}}

\newcommand{\germanA}{{\mathfrak{A}}}
\newcommand{\germanB}{{\mathfrak{B}}}
\newcommand{\germanQ}{{\mathfrak{Q}}}
\newcommand{\germanh}{{\mathfrak{h}}}

\newcommand{\myle}{\preceq}
\newcommand{\myge}{\succeq}
\newcommand{\mygt}{\succ}

\newcommand{\B}{{\sf B}}
\newcommand{\OB}{B^{\rm ord}}
\newcommand{\OS}{{\sf OS}}
\newcommand{\OO}{{\sf O}}
\newcommand{\SP}{{\sf SP}}
\newcommand{\OSP}{{\sf OSP}}
\newcommand{\Eu}{{\sf Eu}}
\newcommand{\ERR}{{\sf ERR}}
\newcommand{\sfB}{{\sf B}}
\newcommand{\sfD}{{\sf D}}
\newcommand{\sfE}{{\sf E}}
\newcommand{\sfG}{{\sf G}}
\newcommand{\sfJ}{{\sf J}}
\newcommand{\sfP}{{\sf P}}
\newcommand{\sfQ}{{\sf Q}}
\newcommand{\sfS}{{\sf S}}
\newcommand{\sfT}{{\sf T}}
\newcommand{\sfW}{{\sf W}}
\newcommand{\sfMV}{{\sf MV}}
\newcommand{\AMV}{{\sf AMV}}
\newcommand{\BM}{{\sf BM}}
\newcommand{\NC}{{\sf NC}}

\newcommand{\emIB}{B^{\rm irr}}
\newcommand{\emIP}{P^{\rm irr}}
\newcommand{\emOB}{B^{\rm ord}}
\newcommand{\emCB}{B^{\rm cyc}}
\newcommand{\emSC}{P^{\rm cyc}}

\newcommand{\stat}{{\rm stat}}
\newcommand{\cyc}{{\rm cyc}}
\newcommand{\Asc}{{\rm Asc}}
\newcommand{\asc}{{\rm asc}}
\newcommand{\Des}{{\rm Des}}
\newcommand{\des}{{\rm des}}
\newcommand{\Exc}{{\rm Exc}}
\newcommand{\exc}{{\rm exc}}
\newcommand{\Wex}{{\rm Wex}}
\newcommand{\wex}{{\rm wex}}
\newcommand{\Fix}{{\rm Fix}}
\newcommand{\fix}{{\rm fix}}
\newcommand{\lrmax}{{\rm lrmax}}
\newcommand{\rlmax}{{\rm rlmax}}
\newcommand{\Rec}{{\rm Rec}}
\newcommand{\rec}{{\rm rec}}
\newcommand{\Arec}{{\rm Arec}}
\newcommand{\arec}{{\rm arec}}
\newcommand{\ERec}{{\rm ERec}}
\newcommand{\erec}{{\rm erec}}
\newcommand{\EArec}{{\rm EArec}}
\newcommand{\earec}{{\rm earec}}
\newcommand{\recarec}{{\rm recarec}}
\newcommand{\nonrec}{{\rm nonrec}}
\newcommand{\Cpeak}{{\rm Cpeak}}
\newcommand{\cpeak}{{\rm cpeak}}
\newcommand{\Cval}{{\rm Cval}}
\newcommand{\cval}{{\rm cval}}
\newcommand{\Cdasc}{{\rm Cdasc}}
\newcommand{\cdasc}{{\rm cdasc}}
\newcommand{\Cddes}{{\rm Cddes}}
\newcommand{\cddes}{{\rm cddes}}
\newcommand{\cdrise}{{\rm cdrise}}
\newcommand{\cdfall}{{\rm cdfall}}
\newcommand{\Peak}{{\rm Peak}}
\newcommand{\peak}{{\rm peak}}
\newcommand{\Val}{{\rm Val}}
\newcommand{\val}{{\rm val}}
\newcommand{\Dasc}{{\rm Dasc}}
\newcommand{\dasc}{{\rm dasc}}
\newcommand{\Ddes}{{\rm Ddes}}
\newcommand{\ddes}{{\rm ddes}}
\newcommand{\inv}{{\rm inv}}
\newcommand{\maj}{{\rm maj}}
\newcommand{\rs}{{\rm rs}}
\newcommand{\cross}{{\rm cr}}
\newcommand{\crosshat}{{\widehat{\rm cr}}}
\newcommand{\nest}{{\rm ne}}
\newcommand{\rodd}{{\rm rodd}}
\newcommand{\reven}{{\rm reven}}
\newcommand{\lodd}{{\rm lodd}}
\newcommand{\leven}{{\rm leven}}
\newcommand{\sg}{{\rm sg}}
\newcommand{\bl}{{\rm bl}}
\newcommand{\tran}{{\rm tr}}
\newcommand{\area}{{\rm area}}
\newcommand{\ret}{{\rm ret}}
\newcommand{\peaks}{{\rm peaks}}
\newcommand{\hl}{{\rm hl}}
\newcommand{\sll}{{\rm sl}}
\newcommand{\negg}{{\rm neg}}
\newcommand{\imp}{{\rm imp}}
\newcommand{\osg}{{\rm osg}}
\newcommand{\ons}{{\rm ons}}
\newcommand{\isg}{{\rm isg}}
\newcommand{\ins}{{\rm ins}}
\newcommand{\LL}{{\rm LL}}
\newcommand{\height}{{\rm ht}}
\newcommand{\as}{{\rm as}}

\newcommand{\ba}{{\bm{a}}}
\newcommand{\bahat}{{\widehat{\bm{a}}}}
\newcommand{\sfa}{{{\sf a}}}
\newcommand{\bb}{{\bm{b}}}
\newcommand{\bc}{{\bm{c}}}
\newcommand{\bchat}{{\widehat{\bm{c}}}}
\newcommand{\bd}{{\bm{d}}}
\newcommand{\bee}{{\bm{e}}}
\newcommand{\beh}{{\bm{eh}}}
\newcommand{\bff}{{\bm{f}}}
\newcommand{\bg}{{\bm{g}}}
\newcommand{\bh}{{\bm{h}}}
\newcommand{\bll}{{\bm{\ell}}}
\newcommand{\bp}{{\bm{p}}}
\newcommand{\br}{{\bm{r}}}
\newcommand{\bs}{{\bm{s}}}
\newcommand{\bt}{{\bm{t}}}
\newcommand{\bu}{{\bm{u}}}
\newcommand{\bw}{{\bm{w}}}
\newcommand{\bx}{{\bm{x}}}
\newcommand{\by}{{\bm{y}}}
\newcommand{\bz}{{\bm{z}}}
\newcommand{\bA}{{\bm{A}}}
\newcommand{\bB}{{\bm{B}}}
\newcommand{\bC}{{\bm{C}}}
\newcommand{\bE}{{\bm{E}}}
\newcommand{\bF}{{\bm{F}}}
\newcommand{\bG}{{\bm{G}}}
\newcommand{\bH}{{\bm{H}}}
\newcommand{\bI}{{\bm{I}}}
\newcommand{\bJ}{{\bm{J}}}
\newcommand{\bM}{{\bm{M}}}
\newcommand{\bN}{{\bm{N}}}
\newcommand{\bP}{{\bm{P}}}
\newcommand{\bQ}{{\bm{Q}}}
\newcommand{\bR}{{\bm{R}}}
\newcommand{\bS}{{\bm{S}}}
\newcommand{\bT}{{\bm{T}}}
\newcommand{\bW}{{\bm{W}}}
\newcommand{\bX}{{\bm{X}}}
\newcommand{\bY}{{\bm{Y}}}
\newcommand{\bIB}{{\bm{B}^{\rm irr}}}
\newcommand{\bOB}{{\bm{B}^{\rm ord}}}
\newcommand{\bOS}{{\bm{OS}}}
\newcommand{\bERR}{{\bm{ERR}}}
\newcommand{\bSP}{{\bm{SP}}}
\newcommand{\bMV}{{\bm{MV}}}
\newcommand{\bBM}{{\bm{BM}}}
\newcommand{\balpha}{{\bm{\alpha}}}
\newcommand{\bbeta}{{\bm{\beta}}}
\newcommand{\bgamma}{{\bm{\gamma}}}
\newcommand{\bdelta}{{\bm{\delta}}}
\newcommand{\bkappa}{{\bm{\kappa}}}
\newcommand{\bmu}{{\bm{\mu}}}
\newcommand{\bomega}{{\bm{\omega}}}
\newcommand{\bphi}{{\bm{\phi}}}
\newcommand{\bsigma}{{\bm{\sigma}}}
\newcommand{\btau}{{\bm{\tau}}}
\newcommand{\bpsi}{{\bm{\psi}}}
\newcommand{\bzeta}{{\bm{\zeta}}}
\newcommand{\bone}{{\bm{1}}}
\newcommand{\bzero}{{\bm{0}}}

\newcommand{\Cbar}{{\overline{C}}}
\newcommand{\Dbar}{{\overline{D}}}
\newcommand{\dbar}{{\overline{d}}}
\def\Ctilde{{\widetilde{C}}}
\def\Ftilde{{\widetilde{F}}}
\def\Gtilde{{\widetilde{G}}}
\def\Htilde{{\widetilde{H}}}
\def\Ptilde{{\widetilde{P}}}
\def\Chat{{\widehat{C}}}
\def\ctilde{{\widetilde{c}}}
\def\zbar{{\overline{Z}}}
\def\pitilde{{\widetilde{\pi}}}

\newcommand{\sech}{{\rm sech}}

%
%
\newcommand{\sn}{{\rm sn}}
\newcommand{\cn}{{\rm cn}}
\newcommand{\dn}{{\rm dn}}
\newcommand{\sm}{{\rm sm}}
\newcommand{\cm}{{\rm cm}}

%
%
\newcommand{\zfz}{ {{}_0 \! F_0} }
\newcommand{\zfo}{ {{}_0  F_1} }
\newcommand{\ofz}{ {{}_1 \! F_0} }
\newcommand{\ofo}{ {{}_1 \! F_1} }
\newcommand{\oft}{ {{}_1 \! F_2} }

%
%
\newcommand{\FHyper}[2]{ {\tensor[_{#1 \!}]{F}{_{#2}}\!} }
\newcommand{\FHYPER}[5]{ {\FHyper{#1}{#2} \!\biggl(
   \!\!\begin{array}{c} #3 \\[1mm] #4 \end{array}\! \bigg|\, #5 \! \biggr)} }
\newcommand{\tfo}{ {\FHyper{2}{1}} }
\newcommand{\tfz}{ {\FHyper{2}{0}} }
\newcommand{\threefz}{ {\FHyper{3}{0}} }
\newcommand{\FHYPERbottomzero}[3]{ {\FHyper{#1}{0} \hspace*{-0mm}\biggl(
   \!\!\begin{array}{c} #2 \\[1mm] \hbox{---} \end{array}\! \bigg|\, #3 \! \biggr)} }
\newcommand{\FHYPERtopzero}[3]{ {\FHyper{0}{#1} \hspace*{-0mm}\biggl(
   \!\!\begin{array}{c} \hbox{---} \\[1mm] #2 \end{array}\! \bigg|\, #3 \! \biggr)} }

\newcommand{\phiHyper}[2]{ {\tensor[_{#1}]{\phi}{_{#2}}} }
\newcommand{\psiHyper}[2]{ {\tensor[_{#1}]{\psi}{_{#2}}} }
\newcommand{\PhiHyper}[2]{ {\tensor[_{#1}]{\Phi}{_{#2}}} }
\newcommand{\PsiHyper}[2]{ {\tensor[_{#1}]{\Psi}{_{#2}}} }
\newcommand{\phiHYPER}[6]{ {\phiHyper{#1}{#2} \!\left(
   \!\!\begin{array}{c} #3 \\ #4 \end{array}\! ;\, #5, \, #6 \! \right)\!} }
\newcommand{\psiHYPER}[6]{ {\psiHyper{#1}{#2} \!\left(
   \!\!\begin{array}{c} #3 \\ #4 \end{array}\! ;\, #5, \, #6 \! \right)} }
\newcommand{\PhiHYPER}[5]{ {\PhiHyper{#1}{#2} \!\left(
   \!\!\begin{array}{c} #3 \\ #4 \end{array}\! ;\, #5 \! \right)\!} }
\newcommand{\PsiHYPER}[5]{ {\PsiHyper{#1}{#2} \!\left(
   \!\!\begin{array}{c} #3 \\ #4 \end{array}\! ;\, #5 \! \right)\!} }
\newcommand{\zerophizero}{ {\phiHyper{0}{0}} }
\newcommand{\ophizero}{ {\phiHyper{1}{0}} }
\newcommand{\zphio}{ {\phiHyper{0}{1}} }
\newcommand{\ophio}{ {\phiHyper{1}{1}} }
\newcommand{\tphio}{ {\phiHyper{2}{1}} }
\newcommand{\tphiz}{ {\phiHyper{2}{0}} }
\newcommand{\tPhio}{ {\PhiHyper{2}{1}} }
\newcommand{\opsio}{ {\psiHyper{1}{1}} }

%
%
\newcommand{\stirlingsubset}[2]{\genfrac{\{}{\}}{0pt}{}{#1}{#2}}
\newcommand{\stirlingcycleold}[2]{\genfrac{[}{]}{0pt}{}{#1}{#2}}
\newcommand{\stirlingcycle}[2]{\left[\! \stirlingcycleold{#1}{#2} \!\right]}
\newcommand{\assocstirlingsubset}[3]{{\genfrac{\{}{\}}{0pt}{}{#1}{#2}}_{\! \ge #3}}
\newcommand{\genstirlingsubset}[4]{{\genfrac{\{}{\}}{0pt}{}{#1}{#2}}_{\! #3,#4}}
\newcommand{\irredstirlingsubset}[2]{{\genfrac{\{}{\}}{0pt}{}{#1}{#2}}^{\!\rm irr}}
\newcommand{\euler}[2]{\genfrac{\langle}{\rangle}{0pt}{}{#1}{#2}}
\newcommand{\eulergen}[3]{{\genfrac{\langle}{\rangle}{0pt}{}{#1}{#2}}_{\! #3}}
\newcommand{\eulersecond}[2]{\left\langle\!\! \euler{#1}{#2} \!\!\right\rangle}
\newcommand{\eulersecondgen}[3]{{\left\langle\!\! \euler{#1}{#2} \!\!\right\rangle}_{\! #3}}
\newcommand{\binomvert}[2]{\genfrac{\vert}{\vert}{0pt}{}{#1}{#2}}
\newcommand{\binomsquare}[2]{\genfrac{[}{]}{0pt}{}{#1}{#2}}


\newenvironment{sarray}{
             \textfont0=\scriptfont0
             \scriptfont0=\scriptscriptfont0
             \textfont1=\scriptfont1
             \scriptfont1=\scriptscriptfont1
             \textfont2=\scriptfont2
             \scriptfont2=\scriptscriptfont2
             \textfont3=\scriptfont3
             \scriptfont3=\scriptscriptfont3
           \renewcommand{\arraystretch}{0.7}
           \begin{array}{l}}{\end{array}}

\newenvironment{scarray}{
             \textfont0=\scriptfont0
             \scriptfont0=\scriptscriptfont0
             \textfont1=\scriptfont1
             \scriptfont1=\scriptscriptfont1
             \textfont2=\scriptfont2
             \scriptfont2=\scriptscriptfont2
             \textfont3=\scriptfont3
             \scriptfont3=\scriptscriptfont3
           \renewcommand{\arraystretch}{0.7}
           \begin{array}{c}}{\end{array}}


\newcommand*\circled[1]{\tikz[baseline=(char.base)]{
  \node[shape=circle,draw,inner sep=1pt] (char) {#1};}}
\newcommand{\ostar}{{\circledast}}
\newcommand{\ostarN}{{\,\circledast_{\vphantom{\dot{N}}N}\,}}
\newcommand{\ostarPsi}{{\,\circledast_{\vphantom{\dot{\Psi}}\Psi}\,}}
\newcommand{\starN}{{\,\ast_{\vphantom{\dot{N}}N}\,}}
\newcommand{\starpsi}{{\,\ast_{\vphantom{\dot{\bpsi}}\!\bpsi}\,}}
\newcommand{\starone}{{\,\ast_{\vphantom{\dot{1}}1}\,}}
\newcommand{\startwo}{{\,\ast_{\vphantom{\dot{2}}2}\,}}
\newcommand{\starinfty}{{\,\ast_{\vphantom{\dot{\infty}}\infty}\,}}
\newcommand{\starT}{{\,\ast_{\vphantom{\dot{T}}T}\,}}

\newcommand*{\Scale}[2][4]{\scalebox{#1}{$#2$}}

\newcommand*{\Scaletext}[2][4]{\scalebox{#1}{#2}} 

\clearpage

It is well known
that the set $\scrt_{n+1}$ of rooted trees on the vertex set
$[n+1] \eqdef \{1,\ldots,n+1\}$ has cardinality $(n+1)^n$;
and from the binomial theorem we have the obvious identity
\be
   (n+1)^n  \;=\;  \sum_{k=0}^n \binom{n}{k} \, n^{n-k}
   \;.
\ee
So it is natural to seek a combinatorial explanation of this identity:
Can we find a partition
of $\scrt_{n+1}$ into subsets $\scrt_{n+1,k}$ ($0 \le k \le n$)
such that $|\scrt_{n+1,k}| = \binom{n}{k} \, n^{n-k}$?

A solution to this problem was found two decades ago
by Chauve, Dulucq and Guibert \cite{Chauve_99,Chauve_00}:
they showed that the number of rooted trees on the vertex set $[n+1]$
in~which exactly $k$ children of the root are lower-numbered than the root is
$\binom{n}{k} \, n^{n-k}$ \cite[A071207]{OEIS}.
Their proof was bijective but rather complicated.\footnote{
   In \cite[Section~3]{Chauve_99},
   the same authors also gave a simple algebraic proof
   of the special case $k=0$,
   based on exponential generating functions
   and the Lagrange inversion formula.
}
Here I would like to give a simpler proof.

Let $T(n;i,k,\ell,m)$ be the number of rooted trees on the vertex set $[n+1]$
in which the root is $i$,
the root has $k$ children $< i$ and $\ell$ children $> i$,
and the forest whose roots are the children $< i$ (resp.\ the children $> i$)
has $m$ (resp.\ $n-m$) vertices.
We~can obtain an explicit formula for $T(n;i,k,\ell,m)$ as follows:
Given $i \in [n+1]$, we choose the $k$ children $< i$ in $\binom{i-1}{k}$ ways,
and the $\ell$ children $> i$ in $\binom{n+1-i}{\ell}$ ways.
Then we choose $m-k$ additional vertices for the first forest
from the remaining $n-k-\ell$ vertices, in $\binom{n-k-\ell}{m-k}$ ways.
This also fixes the $n-m-\ell$ additional vertices for the second forest.
And finally, we recall \cite[Proposition~5.3.2]{Stanley_99}
that the number of forests on $m$ total vertices
with $k$ fixed roots is
\be
   \phi_{m,k}
   \;=\;
   \begin{cases}
      1               & \textrm{if $m=k=0$} \\[0.5mm]
      k \, m^{m-k-1}  & \textrm{if $m \ge 1$ and $0 \le k \le m$} \\[0.5mm]
      0               & \textrm{if $k > m$}
   \end{cases}
 \label{def.phi}
\ee
\noindent
\!\!\cite[A232006]{OEIS}.
In the same way, the number of forests on $n-m$ total vertices
with $\ell$ fixed roots is $\phi_{n-m,\ell}$.
It follows that
\be
   T(n;i,k,\ell,m)
   \;=\;
   \binom{i-1}{k} \, \binom{n+1-i}{\ell} \, \binom{n-k-\ell}{m-k}
                                         \, \phi_{m,k} \, \phi_{n-m,\ell}
   \;.
 \label{eq.T}
\ee
This is defined for $n \ge 0$, $1 \le i \le n+1$, $0 \le k \le n$,
$0 \le \ell \le n-k$ and $k \le m \le n-\ell$.
For $n=0$ the only combinatorially feasible parameters are
$i=1$ and $k = \ell = m = 0$,
and in this case we have $T(0;1,0,0,0) = 1$;
so we can assume henceforth that $n \ge 1$.

We now proceed to sum \reff{eq.T} over~$i$ and $m$.
Note that $i$ appears only in the first two factors
on the right-hand side of \reff{eq.T},
while $m$ appears only in the final three factors.
So we can perform these two sums separately.

\bigskip

{\bf Sum over $\bm{i}$.}
We claim that for any integers $n,k,\ell \ge 0$, we have
\be
   \sum_{i=1}^{n+1} \binom{i-1}{k} \, \binom{n+1-i}{\ell}
   \;=\;
   \binom{n+1}{k+\ell+1}
   \;.
 \label{eq.sum1}
\ee
This identity has a simple combinatorial proof:
the right-hand side is the number of ways of choosing
$k+\ell+1$ elements from the set $[n+1]$;
if we arrange these elements in increasing order
and call the $(k+1)$st of them $i$,
then the two binomial coefficients on the left-hand side
give the number of ways of choosing the first $k$ elements
and the last $\ell$ elements, respectively.
The identity \reff{eq.sum1} can also be derived algebraically
as a corollary of the Chu--Vandermonde identity;
we discuss this in Appendix~\ref{app.chu-vandermonde}.

{}From the right-hand side, we see in particular that
\reff{eq.sum1} depends on $k$ and $\ell$ only via their sum.


\bigskip

{\bf Sum over $\bm{m}$.}
We claim that for any integers $n,k,\ell \ge 0$ with $k+\ell \le n$, we have
\be
   \sum_{m=k}^{n-\ell} \binom{n-k-\ell}{m-k} \, \phi_{m,k}
                                             \, \phi_{n-m,\ell}
   \;=\;
   \phi_{n,k+\ell}
   \;.
 \label{eq.sum2}
\ee
This identity too has a simple combinatorial proof:
the right-hand side counts the forests on the vertex set $[n]$
with $k+\ell$ fixed roots,
while the left-hand side partitions this count according to the number $m$
of vertices that belong to the subforest associated to the first $k$ roots.
The identity \reff{eq.sum2} can also be derived algebraically
as a corollary of an Abel identity;
we discuss this in Appendix~\ref{app.abel}.

{}From the right-hand side, we see in particular that
\reff{eq.sum2} depends on $k$ and $\ell$ only via their sum.

\bigskip

{\bf Combining the two sums.}
Combining \reff{eq.T} with \reff{eq.sum1} and \reff{eq.sum2}, we have
for $n \ge 1$
\begin{eqnarray}
   \sum_{i=1}^{n+1}  T(n;i,k,\ell,m)
   & = &
   \binom{n+1}{k+\ell+1} \, \binom{n-k-\ell}{m-k}
                         \, \phi_{m,k} \, \phi_{n-m,\ell}
       \label{eq.ans1}  \\[4mm]
   \sum_{m=k}^{n-\ell}   T(n;i,k,\ell,m)
   & = &
   \binom{i-1}{k} \, \binom{n+1-i}{\ell} \, (k+\ell) \, n^{n-k-\ell-1}
       \label{eq.ans2}  \\[4mm]
   \sum_{i=1}^{n+1}  \sum_{m=k}^{n-\ell}   T(n;i,k,\ell,m)
   & = &
   \binom{n+1}{k+\ell+1} \, (k+\ell) \, n^{n-k-\ell-1}
       \label{eq.ans3}
\end{eqnarray}
The right-hand side of \reff{eq.ans3}
depends on $k$ and $\ell$ only via their sum;
we denote this quantity by $g_n(k+\ell)$, i.e.\ we define
\be
   g_n(K)  \;\eqdef\; \binom{n+1}{K+1} \, K \, n^{n-K-1}
      \quad\hbox{for $n \ge 1$ and $0 \le K \le n$}
   \;.
 \label{def.gnK}
\ee

{\bf Sum over $\bm{\ell}$.}
The final step is to sum \reff{eq.ans3} over $\ell$ at fixed $k$,
i.e.\ to compute
\be
   G_n(k)  \;\eqdef\; \sum_{\ell=0}^{n-k} g_n(k+\ell)
           \;=\;  \sum_{K=k}^n g_n(K)
   \;.
 \label{def.GnK}
\ee
We prove that $G_n(k) = \binom{n}{k} n^{n-k}$, as follows:
{}From \reff{def.GnK}, $G_n(k)$ manifestly satisfies the backward recurrence
\be
   G_n(k)  \;=\; G_n(k+1)  \:+\: \binom{n+1}{k+1} \, k \, n^{n-k-1}
\ee
with initial condition $G_n(n) = 1$.
A simple calculation shows that $\widehat{G}_n(k) = \binom{n}{k} n^{n-k}$
satisfies the same recurrence and the same initial condition.
Hence $G_n(k) = \widehat{G}_n(k)$.
QED

Xi Chen (private communication) has found an alternate proof of
$G_n(k) = \binom{n}{k} n^{n-k}$ that {\em derives}\/ it
(rather than simply pulling it out of a hat, as the foregoing proof does);
this proof is presented in Appendix~\ref{app.chen}.

\bigskip

{\bf Three final remarks.}

1. The special case $k=0$ of \reff{eq.ans3}
was found by Chauve {\em et al.}\/ \cite[Proposition~2]{Chauve_00}.

2. By summing \reff{eq.ans2} over $\ell$,
we can compute the number of rooted trees in $\scrt_{n+1,k}$
that have a specified element $i$ as the root.
This sum is easily performed using the binomial theorem and its derivative,
and gives
\be
   \sum_{\ell=0}^{n+1-i} \sum_{m=k}^{n-\ell}   T(n;i,k,\ell,m)
   \;=\;
   \binom{i-1}{k} \, \big[ (k+1)(n+1) - i \big] \, n^{i-k-2} \, (n+1)^{n-i}
   \;.
\ee
%
%
%
%
%
%
%
For the special case $k=0$, this result was obtained bijectively by
Chauve {\em et al.}\/ \cite[proof of Proposition~1]{Chauve_00}.

3. We can also compute the number of rooted trees on $n+1$ labeled vertices
in which the root has exactly $K$ children:  it suffices to sum
\reff{eq.ans3} over $k,\ell \ge 0$ with $k+\ell = K$, yielding
\be
   (K+1) \, \binom{n+1}{K+1} \, K \, n^{n-K-1}
   \;=\;
   (n+1) \, \binom{n}{K} \, K \, n^{n-K-1}
   \;.
\ee
Here $n+1$ counts the number of choices for the root,
and the remaining factor
$f_{n,k} = \binom{n}{K} \, K \, n^{n-K-1} = \binom{n}{K} \, \phi_{n,K}$
counts the number of $K$-component forests of rooted trees
on $n$~labeled vertices.
This latter result is essentially equivalent to \reff{def.phi},
and is well known.\footnote{
    See e.g.\ \cite{Clarke_58}, \cite[pp.~26--27]{Moon_70},
    \cite[p.~70]{Comtet_74}, \cite[pp.~25--28]{Stanley_99}
    or \cite{Avron_16}.
    See also \cite{Riordan_68b,Shor_95,Zeng_99,Guo_17}
    and \cite[pp.~235--240]{Aigner_18}
    for related information.
}

%

\bigskip

{\bf Note Added:}
After my posting of the preprint version of this manuscript,
Jiang Zeng kindly showed me the following quick and elegant proof
of \reff{eq.ans3}:

We can construct rooted trees on the vertex set $[n+1]$
with $k$ (resp.\ $\ell$) children smaller (resp.\ larger) than the root,
as follows:
Choose a subset $S \subseteq [n+1]$ of cardinality $k+\ell+1$
--- let us call its elements $a_1 < \ldots < a_{k+\ell+1}$ ---
and then construct a tree with $a_{k+1}$ as the root
and the $k+\ell$ elements of $S \setminus \{a_{k+1}\}$ as children of the root.
By \reff{def.phi} there are
\be
   \binom{n+1}{k+\ell+1} \, \phi_{n,k+\ell}
   \;=\;
   \binom{n+1}{k+\ell+1} \, (k+\ell) \, n^{n-k-\ell-1}
\ee
such trees; this is \reff{eq.ans3}.

The longer proof given in the body of this paper
may nevertheless still be of some interest,
as it yields the more refined enumerations \reff{eq.ans1} and \reff{eq.ans2}.

\appendix
\section*{Appendix: Algebraic proofs}
\setcounter{section}{1}
\setcounter{equation}{0}
\def\theequation{\Alph{section}.\arabic{equation}}

\subsection{A corollary of the Chu--Vandermonde identity}
   \label{app.chu-vandermonde}

The identity \reff{eq.sum1} is a special case
of a slightly more general binomial identity, namely
\be
   \sum_{j=k-m}^{n-\ell} \binom{m+j}{k} \binom{n-j}{\ell}
   \;=\;
   \binom{m+n+1}{k+\ell+1}
   \;,
 \label{eq.binomial_identity}
\ee
valid for integers $k,\ell,m,n$ with $k,\ell \ge 0$ and $m+n \ge -1$.
Although this identity can be found
in several places in the literature\footnote{
   See e.g.\ 
   \cite[p.~22, eq.~(3.3)]{Gould_72}
   and
   \cite[p.~169, eq.~(5.26) and pp.~243, 527, Exercise 5.14]{Graham_94}.
},
I have been unable to find any place where it is stated clearly
with its optimal conditions of validity.
I~will therefore give here a detailed derivation,
keeping careful track of the conditions of validity for each step.

The binomial coefficients are defined as usual by \cite[p.~154]{Graham_94}
\be
   \binom{r}{k}
   \;=\;
   \begin{cases}
        \displaystyle {r(r-1) \,\cdots\, (r-k+1) \over k!}
                             & \textrm{for integer $k \ge 0$}
              \\[1mm]
        0                    & \textrm{for integer $k < 0$}
              \\[1mm]
        \textrm{undefined}   &   \textrm{if $k$ is not an integer}
   \end{cases}
\ee
Here $r$ can be any element of any commutative ring containing the rationals;
in particular, it can be an indeterminate in a ring of polynomials
over the rationals.
The binomial coefficients satisfy
\be
   \binom{r}{k}
   \;=\;
   (-1)^k \, \binom{-(r-k+1)}{k}
      \qquad\hbox{for integer $k$}
 \label{eq.upper_negation}
\ee
(``upper negation'') and
\be
   \binom{n}{k}  \;=\;  \binom{n}{n-k}
      \qquad\hbox{for integer $n \ge 0$ and integer $k$}
 \label{eq.symmetry}
\ee
(``symmetry'').
Finally, they satisfy the {\em Chu--Vandermonde identity}\/
\be
   \sum_{j=0}^N \binom{x}{j} \binom{y}{N-j}  \;=\;  \binom{x+y}{N}
      \qquad\hbox{for integer $N$}
   \;,
 \label{eq.vandermonde}
\ee
where $x$ and $y$ can be indeterminates.
Applying \reff{eq.upper_negation} to all three binomial coefficients
in the Chu--Vandermonde identity
and then replacing $x \to -x$, $y \to -y$,
we obtain the {\em dual Chu--Vandermonde identity}\/
\be
   \sum_{j=0}^N \binom{x+j-1}{j} \binom{y+N-j-1}{N-j}
   \;=\;
   \binom{x+y+N-1}{N}
      \qquad\hbox{for integer $N$}
   \;.
 \label{eq.dual_vandermonde}
\ee

Now suppose that $x,y$ are integers $\ge 1$ and that $x+y+N \ge 1$;
then we can apply the symmetry \reff{eq.symmetry}
to the three binomial coefficients in \reff{eq.dual_vandermonde}.
Writing $x=k+1$ and $y = \ell+1$ with integers $k,\ell \ge 0$, we have
\begin{eqnarray}
   & &
   \sum_{j=0}^N \binom{k+j}{k} \binom{N+\ell-j}{\ell}
   \;=\;
   \binom{k+\ell+N+1}{k+\ell+1}
        \nonumber \\[1mm]
   & &
      \qquad\hbox{for integers $k,\ell,N$ with $k,\ell \ge 0$
                          and $k+\ell+N \ge -1$}
   \;.
\end{eqnarray}
Now change variables $j = j'+m-k$ and $N = m+n-k-\ell$:
\begin{eqnarray}
   & &
   \sum_{j'=k-m}^{n-\ell} \binom{m+j'}{k} \binom{n-j'}{\ell}
   \;=\;
   \binom{m+n+1}{k+\ell+1}
        \nonumber \\[1mm]
   & &
      \qquad\hbox{for integers $k,\ell,m,n$ with $k,\ell \ge 0$
                          and $m+n \ge -1$}
   \;.
\end{eqnarray}
Dropping primes, this is \reff{eq.binomial_identity}.


\subsection{Abel identity}   \label{app.abel}

The identity \reff{eq.sum2} can also be derived algebraically, as follows:
We begin from the well-known Abel identity \cite[p.~73]{Roman_84}
\be
   \sum_{M=0}^N \binom{N}{M} \, x (x+M)^{M-1} \, y (y+N-M)^{N-M-1}
   \;=\;
   (x+y) \, (x+y+N)^{N-1}
 \label{eq.abel.2}
\ee
(see also \cite[p.~20, eq.~(20)]{Riordan_68} multiplied by $xy$).\footnote{
   The identity \reff{eq.abel.2} asserts that
   the polynomials $P_N(x) = x (x+N)^{N-1}$,
   which are a specialization of the celebrated {\em Abel polynomials}\/
   $A_n(x;a) = x (x-an)^{n-1}$ \cite{Mullin_70,Francon_74,Sagan_83,Roman_84}
   to $a=-1$,
   form a {\em sequence of binomial type}\/
   \cite{Mullin_70,Garsia_73,Roman_84}.
   See also \cite{Labelle_81} \cite[Section~3.1]{Bergeron_98}
   for a purely combinatorial approach to sequences of binomial type,
   employing the theory of species.
}
Since all the terms in this identity (even the ones with $M=0$ and $M=N$)
are polynomials in $x$ and~$y$,
the variables $x$ and $y$ can be specialized without restriction.
(Note, however, that in applying this identity,
we must first fix $N$ and $M$ and then specialize $x$ and $y$.)
Setting $N = n-k-\ell$ and changing variables by $M = m-k$ yields
\begin{eqnarray}
   & &
   \sum_{m=k}^{n-\ell} \binom{n-k-\ell}{m-k} \, x (x+m-k)^{m-k-1}
                                             \, y (y+n-m-\ell)^{n-m-\ell-1}
   \qquad\qquad
          \nonumber \\[2mm]
   & & \hspace*{4cm}
   \;=\;
   (x+y) \, (x+y+n-k-\ell)^{n-k-\ell-1}
   \;.
 \label{eq.sum2pre}
\end{eqnarray}
Specializing now to $x=k$ and $y=\ell$,
we see that $\left. x (x+m-k)^{m-k-1} \right|_{x=k} = \phi_{m,k}$
even when $m=k=0$,
and likewise
$\left. y (y+n-m-\ell)^{n-m-\ell-1} \right|_{y=\ell} = \phi_{n-m,\ell}$
even when $n-m=\ell=0$.
It follows that
\be
   \sum_{m=k}^{n-\ell} \binom{n-k-\ell}{m-k} \, \phi_{m,k}
                                             \, \phi_{n-m,\ell}
   \;=\;
   (k+\ell) \, n^{n-k-\ell-1}
   \;=\;
   \phi_{n,k+\ell}
   \;,
 \label{eq.sum2bis}
\ee
valid for $n \ge 1$ and $k,\ell \ge 0$ with $k+\ell \le n$.

We remark, finally, that many Abel identities, including \reff{eq.abel.2},
can be proven combinatorially:
see e.g.\ \cite{Francon_74,Sagan_83,Pitman_02}.

\subsection{Alternate proof of $G_n(k) = \binom{n}{k} n^{n-k}$
   (due to Xi~Chen)}   \label{app.chen}

We compute the row-generating polynomials
$\scrg_n(x) \eqdef \sum\limits_{k=0}^n G_n(k) \, x^k$, as follows:
\begin{subeqnarray}
   \scrg_n(x)
   & = &
   \sum_{k=0}^n \sum_{K=k}^n \binom{n+1}{K+1} \, K \, n^{n-K-1} \, x^k
       \\[2mm]
   & = &
   n^{n-1}
   \sum_{K=0}^n \binom{n+1}{K+1} \, K \, n^{-K} \sum_{k=0}^K \, x^k
       \\[2mm]
   & = &
   n^{n-1}
   \sum_{K=0}^n \binom{n+1}{K+1} \, K \, n^{-K} \, {1 - x^{K+1} \over 1-x}
       \\[2mm]
   & = &
   {n^{n-1} \over 1-x}
   \left[ \sum_{K=0}^n \binom{n+1}{K+1} \, K \, {1 \over n^K}
          \:-\:
          x \sum_{K=0}^n \binom{n+1}{K+1} \, K \, {x^K \over n^K}
   \right]
       \\[2mm]
   & = &
   {n^{n-1} \over 1-x}
   \big[ \scrf_n(1/n) \:-\: x \, \scrf_n(x/n) \big]
 \label{eqns.scrg}
\end{subeqnarray}
where
\be
   \scrf_n(x)  \;\eqdef\;  \sum_{K=0}^n \binom{n+1}{K+1} \, K \, x^K
   \;.
\ee
A simple computation, using the derivative of the binomial theorem,
shows that
\be
   \scrf_n(x)
   \;=\;
   (n+1) (x+1)^n  \,-\, {(x+1)^{n+1} - 1  \over x}
   \;.
\ee
Therefore
\be
   \scrf_n(1/n) \;=\;  n
   \qquad\hbox{and}\qquad
   x \, \scrf_n(x/n) \;=\;  {1 \over n^{n-1}} (x-1)(x+n)^n \,+\, n
   \;,
\ee
and inserting these into \reff{eqns.scrg} gives
\be
   \scrg_n(x)  \;=\;  (x+n)^n
   \;.
\ee
Taking the coefficient of $x^k$ in $\scrg_n(x)$,
we conclude that $G_n(k) = \binom{n}{k} n^{n-k}$.


\section*{Acknowledgments}

I wish to thank Xi Chen and Alex Dyachenko for helpful conversations,
and Jiang Zeng for correspondence.
I am especially grateful to Xi Chen for giving me permission
to include her proof of $G_n(k) = \binom{n}{k} n^{n-k}$
in Appendix~\ref{app.chen},
and to Jiang Zeng for giving me permission to include his quick proof
of \reff{eq.ans3}.

This work was immeasurably facilitated by
the On-Line Encyclopedia of Integer Sequences \cite{OEIS}.
I warmly thank Neil Sloane for founding this indispensable resource,
and the hundreds of volunteers for helping to maintain and expand it.

This research was supported in part by
Engineering and Physical Sciences Research Council grant EP/N025636/1.

\end{document}